\documentclass[12pt]{article}
\usepackage{amsmath,amssymb}

\newcommand{\N}{{\mathbb N}}
\newcommand{\R}{{\mathbb R}}
\newcommand{\Z}{{\mathbb Z}}

\newcommand{\A}{\mathcal A}

\newcommand{\LL}{{\mathcal L}}

\newcommand{\Cov}{{\rm Cov}}

\newcommand{\cla}{\stackrel{{\cal L}}{\rightarrow}}

\catcode`@=11
\@addtoreset{equation}{section}
\catcode`@=12

\newtheorem{lemma}{Lemma}[section]
\newtheorem{proposition}[lemma]{Proposition}
\newtheorem{corollary}[lemma]{Corollary}

\newtheorem{remark}[lemma]{Remark}
\newtheorem{theorem}{Theorem}

\begin{document}

\title
{New Techniques for Empirical Processes of\\ Dependent Data}
\author{\sc Herold Dehling\thanks{Fakult\"at f\"ur Mathematik, 
Ruhr-Universit\"at Bochum, Universit\"atsstra\ss e 150, 44780 Bochum, Germany; 
e-mail: herold.dehling@rub.de} \and
\sc Olivier Durieu\thanks{Laboratoire de Mathematique Rapha\"el Salem, UMR 6085 CNRS-Universit\'e de Rouen,
e-mail: olivier.durieu@etu.univ-rouen.fr}
\and
\sc Dalibor Volny\thanks{Laboratoire de Mathematique Rapha\"el Salem, UMR 6085 CNRS-Universit\'e de Rouen, e-mail:
dalibor.volny@univ-rouen.fr}}
\maketitle

\begin{abstract}
We present a new technique for proving empirical process invariance principle for stationary processes
$(X_n)_{n\geq 0}$.
The main novelty of our approach lies in the fact that we only require the central limit theorem 
and a moment bound for a restricted class of functions $(f(X_n))_{n\geq 0}$, not containing the indicator 
functions.
Our approach can be applied to Markov chains and dynamical systems, using spectral properties of the transfer operator.
Our proof consists of a novel application of chaining techniques.
\end{abstract}

\section{Introduction}
Let $(X_n)_{n\geq 0}$ be a stationary ergodic process of $\R$-valued
random variables with marginal distribution function $F(t)=P(X_0\leq t)$. Define
the empirical distribution function $(F_n(t))_{t\in\R}$ and the 
empirical process $(U_n(t))_{t\in\R}$ by
\begin{eqnarray*}
 F_n(t) &:= & \frac{1}{n} \sum_{i=1}^n 1_{(-\infty,t]} (X_i),\; t\in\R,
 \\
 U_n(t) &:=& \sqrt{n} (F_n(t)-F(t)),\; t\in\R.
\end{eqnarray*}
The empirical process plays a prominent role in non-parametric statistical 
inference about the distribution function $F$. In all statistical  
applications, information about the distribution of the
empirical process is needed. 

\medskip

In the case of i.i.d. observations, Donsker \cite{Don52} proved in 1952 that  
the empirical process converges in distribution to a Brownian bridge 
process, thus confirming an earlier conjecture of Doob \cite{Doo49}. In 1968,
Billingsley \cite{Bil68} extended Donsker's theorem to
some weakly dependent processes, specifically to
functionals of $\phi$-mixing processes. One of the applications of 
Billingsley's theorem is to the empirical process of data generated by the 
continued fraction 
dynamical system $T:[0,1]\rightarrow [0,1]$, $T(x):=\frac{1}{x}$.
Since 1968, many authors have studied the empirical process of weakly
dependent data. Invariance principles for empirical distribution of strong mixing random variables were proved in 1977
by Berkes and Philipp \cite{BerPhi77} and in 1980 for the multivariate case by Philipp and Pinzur \cite{PhiPin80}.
Later, absolutely regular processes were studied by Doukhan et al. \cite{DouMasRio95} and Borovkova et al. \cite{BorBurDeh01}. Many other weak dependence conditions have been studied, see Doukhan and Louichi \cite{DouLou99},
Prieur \cite{Pri02}, Dedecker and Prieur \cite{DedPri07}, Wu and Shao \cite{WuSha04}, Wu \cite{Wu08}.
From the point of view of dynamical systems, an empirical process invariance principle for expanding maps of the interval was proved by Collet et al \cite{ColMarSch04}. Another one for ergodic torus automorphisms was proved by Durieu and Jouan \cite{DurJou08}.

\medskip

Proofs of empirical process invariance principles usually consist of two parts,
establishing finite-dimensional convergence and tightness of the empirical 
process. Finite-dimensional convergence, i.e. convergence in distribution
of the sequence
of vectors $(U_n(t_1),\ldots,U_n(t_k))_{n\geq 1}$, is  an immediate 
consequence of the multivariate CLT for partial sums of the process
\[
 (1_{(-\infty,t_1]}(X_n),\ldots,1_{(-\infty,t_k]}(X_n))_{n\geq 1}.
\]
Tightness is far more difficult to establish. One ingredient is usually a
probability bound on the increments of the empirical process 
\[
 U_n(t)-U_n(s) = \frac{1}{\sqrt{n}} \sum_{i=1}^n\{ 1_{(s,t]}(X_i)-(F(t)-F(s))\},
\]
for a fixed pair $ s <t $.
Such bounds can in the simplest approach be obtained from bounds on the 4-th 
moments of $U_n(t)-U_n(s)$. Other results require higher order moment bounds
or even exponential bounds. 

\medskip

The traditional approach to empirical process invariance principles, as 
outlined above, works well in situations when the sequence of
indicator variables 
$(1_{(s,t]}(X_n))_{n\geq 0}$ inherits good properties from the original 
process  $(X_n)_{n \geq 0}$. This holds, for example, when $(X_n)_{n\geq 0}$
is strong (uniform, beta) mixing, because then $(1_{(s,t]}(X_n))_{n\geq 0}$ has
the same property. There are, however, situations where this is not the case
or at least not easy to establish. For some types of 
Markov processes and dynamical systems, see e.g. Hennion and Herv\'e \cite{HenHer01}, 
one has good control over the 
properties of $(f(X_n))_{n\geq 0}$ when $f$ is a Lipschitz function, but not
for indicator functions.
For example, Gou\"ezel \cite{Gou08} gave an uniformly expanding map of the interval
which has a spectral gap on the space of Lipschitz functions but not on the space of bounded variation functions.
In this paper, we develop an approach that is strictly based on properties of
Lipschitz functions $f(X_i)$ of the original data. 
We make two basic assumptions, namely that the partial sums of Lipschitz
functions satisfy the CLT and that a suitable 4-th moment bound is 
satisfied. 

\medskip

For our proof we develop a variant of the classical
chaining technique that uses only Lipschitz 
functions at all stages of the chaining argument. We replace the usual
finite-dimensional convergence plus tightness approach by a method
of approximation by a sequence of 
finite-dimensional processes, which are different from
the coordinate projections $(U_n(t_1),\ldots,U_n(t_k))$. We show convergence
in distribution of the finite-dimensional processes and prove that the
finite-dimensional process approximates the empirical process. In the final
step, we use an improved version of a Theorem of Billingsley \cite{Bil68}, see our 
Theorem \ref{pbil} below, to establish convergence in distribution of the 
empirical process.

\medskip

In the present paper, we make two assumptions concerning the process 
$(X_i)_{i\geq 0},$ 

\begin{enumerate}
\item 
For any Lipschitz function $f$, the CLT holds, i.e.
\begin{equation}
 \frac{1}{\sqrt{n}} \sum_{i=1}^{n}\{f(X_i)-Ef(X_i)\}
  \cla N(0,\sigma^2),
\label{eq:lip-clt}
\end{equation}
where $N(0,\sigma^2)$ denotes a normal law with 
mean zero and variance 
\[
\sigma^2=E(f(X_0)-Ef(X_0))^2+2\sum_{i=1}^\infty\Cov(f(X_0),f(X_i)).
\]

\item 
A bound on the 4-th central moments of partial sums of $(f(X_i))_{i\geq0}$, 
$f$ bounded Lipschitz with $E(f(X_0))=0$, of the type
\begin{eqnarray}
&&E\left\{\sum_{i=1}^{n}f(X_i)\right\}^4 \nonumber\\
&&\quad \leq C m_f^3 \left( n \|f(X_0) \|_1 
\log^\alpha\left(1+ \| f \|\right) + n^2 \|f (X_0)\|^2_1 \log^\beta\left(1+ \|f \|\right)\right),\nonumber\\
&&
\label{eq:4th-moment}
\end{eqnarray}
where C is some universal constant, $\alpha$ and $\beta$ are some nonnegative integers,
\begin{equation*}
 \parallel f \parallel 
 = \sup_x \mid f(x) \mid + \sup_{x\neq y} 
 \frac{\mid f(x)- f(y) \mid}{\mid x-y \mid}
\end{equation*}
and
$$
m_f=\max\{1,\sup_x \mid f(x) \mid\}.
$$
\end{enumerate}

\begin{remark}
These assumptions can be verified for a large class of Markov chains and dynamical systems.
Concerning the CLT for Lipschitz functions, many results can be found in the literature; see e.g. Hennion and Herv\'e \cite{HenHer01} for the spectral gap method and Bradley \cite{Bra07} for mixing approach.
Durieu \cite{Dur08b} has established 4-th moment bounds of the type (\ref{eq:4th-moment}) for Markov chains and dynamical systems under spectral properties. For more details and concrete examples see Section \ref{examp} of this paper.
\end{remark}

We shall assume some regularity for the distribution function of $X_0$. We define the modulus of continuity of a function
$f:\R\longrightarrow \R$ by
$$
\omega_f(\delta)=\sup\left\{|f(s)-f(t)|\,:\,s,t\in\R,|s-t|<\delta\right\}.
$$

We can now state our main result.

\begin{theorem}\label{thm1}
Let $(X_i)_{i\ge 0}$ be an $\R$-valued stationary ergodic random process such that
the conditions (\ref{eq:lip-clt}) and (\ref{eq:4th-moment}) hold. Assume that $X_0$ has a distribution function $F$ satisfying the following condition,
\begin{equation}\label{eq-modulus}
 \omega_F(\delta)\leq D|\log(\delta)|^{-\gamma}\mbox{ for some } D>0 \mbox{ and } \gamma>\max\{\frac{\alpha}{2},\beta\},
\end{equation}
then
\begin{equation*}
 (U_n(t))_{t\in\R} \stackrel{\mathcal{D}}
 {\longrightarrow} (W(t))_{t\in\R},
\end{equation*}
where $W(t)$ is a mean-zero Gaussian process with covariances
\begin{eqnarray*} 
 E  W (s) \cdot W (t) 
  &=& \Cov ( 1_{(- \infty , s]}(X_0) , 1_{( - \infty, t]} (X_0))\\
  &&+ \sum^\infty_{k=1} \Cov(1_{(-\infty, s]}(X_0), 1_{(-\infty, t]}(X_k))\\
  &&+ \sum^\infty_{k=1} \Cov (1_{(-\infty, s]}(X_k), 1_{(-\infty, t]}(X_0)).
\end{eqnarray*}
Further, almost surely, $(W(t))_{t\in\R}$ has continuous sample paths.
\end{theorem}

\begin{remark}
 In particular, if $X_0$ has a H\"older-continuous distribution function then (\ref{eq-modulus}) holds.
\end{remark}
\begin{remark}
 If the $X_i$'s are i.i.d., $(W(t))_{t\in\R}$ is a Brownian bridge, but this is not the case for dependent variables,
as in Billingsley \cite{Bil68} or Collet et al. \cite{ColMarSch04}.
\end{remark}

\medskip

In order to prove Theorem \ref{thm1}, we apply the following theorem, which is a stronger version of Theorem 4.2 of 
Billingsley \cite{Bil68} in the complete case. We do not need to assume a priori that $X^{(m)}$ 
has a limit in distribution.

\begin{theorem}\label{pbil}
Let $ (S,\rho)$ be a complete separable metric space and let 
$X_n,X_n^{(m)}$ and $X^{(m)}$, $n,m \geq 1$ be S-valued random variables 
satisfying
\begin{eqnarray}
 && X_n^{(m)}\stackrel {\mathcal{D}}{\longrightarrow } X^{(m)} \mbox{ as }  
 n \rightarrow \infty , \forall m \label{eq:p1}\\ 
 && \lim_{m \rightarrow \infty} \limsup_ {n \rightarrow \infty} 
 P ( \rho (X_n, X_n^{(m)}) \geq \varepsilon) = 0, 
 \forall \varepsilon > 0 .\label{eq:p2}
\end{eqnarray}
Then there exists an S-valued random variable $X$ such that
\[ 
  X_n \stackrel{\mathcal{D}}{\longrightarrow} X \mbox{ as } 
 n \rightarrow \infty .
\]
Moreover  $ X^{(m)} \stackrel{\mathcal{D}} {\longrightarrow} 
 X \mbox{ as } m \rightarrow \infty$.
\end{theorem}

Both theorems are proved in Section 2 and Section 3.

\section{Proof of Theorem \ref{thm1}}

\subsection{The bounded case}\label{boundedcase}

We first prove the result for bounded variables. 
Let $(X_i)_{i\ge 0}$ be a $[0,1]$-valued stationary ergodic random process such that
(\ref{eq:lip-clt}), (\ref{eq:4th-moment}) and (\ref{eq-modulus}) hold.

\medskip

In our approach we work with Lipschitz continuous approximations to the 
indicator functions 
$1_{(-\infty,t]}(x)$. Given a partition 
\begin{equation*}
0=t_0^\prime< \ldots <t_m^\prime=1
\end{equation*}
we define
$$
t_j=F^{-1}(t_j^\prime)
$$
where $F^{-1}$ is given by
$$
F^{-1}(t)=\sup\{s\in [0,1] : F(s)\leq t\}.
$$
Thus, by continuity of $F$, we have a partition
$$
0\le t_0<\dots<t_m=1.
$$
We introduce the functions $\varphi_j:[0,1]\rightarrow\mathbb{R}$
by
\begin{equation*}
\varphi_j(x)=\varphi\left(\frac{x-t_{j-1}}{t_{j-1}-t_{j-2}}\right),\quad \mbox{ for } j=2,\dots,m
\end{equation*}
where 
\begin{equation}
 \varphi(x)=1_{(-\infty,-1]}(x)-x1_{(-1,0]}(x)
\label{eq:phi}
\end{equation} 
and $\varphi_1\equiv0$.\\
The function 
$\varphi_j$ will serve as a Lipschitz-continuous approximation to the 
indicator function $1_{(-\infty,t_{j-1}]}(x).$ Note that 
$\varphi_j (x) $ depends on the partition, not only on the point $t_{j-1}.$
We now define the process
\begin{eqnarray*}
 F_n^{(m)}(t) 
 &=& \frac{1}{n} \sum_{i=1}^{n} 
   \sum_{j=1}^{m} 1_{[t_{j-1}, t_j)}(t) \varphi_j(X_i)\\
 &=& \sum_{j=1}^{m}\left(\frac{1}{n} 
   \sum_{i=1}^{n} \varphi_j(X_i)\right) 1_{[t_{j-1}, t_j)}(t).
\end{eqnarray*}
Note that $F_n^{(m)}(t)$ is a piecewise constant approximation to the 
empirical distribution function $F_n(t)$. For $t\in [t_{j-1},t_j]$,
we have the inequality
\[
  F_n(t_{j-2}) \leq F_n^{(m)} (t) \leq F_n(t_{j-1}).
\]
We define further
\[
  F^{(m)}(t) = E \left(F_n^{(m)}(t)\right) = 
 \sum_{j=1}^m E \left(\varphi_j(X_0)\right) 1_{[t_{j-1},t_j)}(t),
\]
and finally the centered and normalized process
\begin{equation}
  U_n^{(m)}(t) = \sqrt{n} \left(F_n^{(m)}(t) - F^{(m)}(t)\right).
\label{eq:unm} 
\end{equation}
Our proof of Theorem \ref{thm1} now consists of two parts, each of which will be 
formulated separately as a proposition below. The theorem will follow by application of Theorem
\ref{pbil}, where $(S,\rho)$ is the space of cadlag functions $D[0,1]$ provided with the Skorohod topology and the metric
$d_0$; see Billingsley \cite{Bil68}, p. 113. Note that $(D[0,1],d_0)$ is a complete separable metric space.

\begin{proposition} 
For any partition $0=t_0^\prime < \ldots < t_m^\prime=1$, there exists a 
piecewise constant Gaussian process
$\left(W^{(m)}(t)\right)_{0\leq t\leq 1}$ such that
\[
 \left(U_n^{(m)}(t)\right)_{0\leq t \leq 1} \mathop{\longrightarrow}
 \limits^{\cal{D}} \left(W^{(m)}(t)\right)_{0\leq t \leq 1}.
\]
The sample paths of the processes $\left(W^{(m)}(t)\right)_{0 \leq t \leq 1}$ 
are  constant on each of the intervals
$[t_{j-1}, t_j)$, $1 \leq j \leq m, $ and $W^{(m)}(0) = 0.$
The vector $(W^{(m)} (t_1), \ldots, W^{(m)}(t_m))$
has a multivariate normal distribution with mean zero and covariances
\begin {eqnarray*}
\Cov (W^{(m)}(t_{i-1}), W^{(m)}(t_{j-1}))
&=& \Cov (\varphi_i(X_0), \varphi_j (X_0))\\
&&+ \sum^\infty_{k=1}\Cov (\varphi_i(X_0), \varphi_j (X_k))\\
&&+ \sum^\infty_{k=1}\Cov (\varphi_i(X_k), \varphi_j (X_0))
\end{eqnarray*}
\label{prop:fidi-conv}
\end{proposition}

{\em Proof.}

Using (\ref{eq:lip-clt}) and the Cram\'er-Wold device, we can show that
for any Lipschitz functions $f_1, \ldots ,f_k,$ the multivariate 
CLT holds, i.e. 
\begin{equation*}
 \frac{1}{\sqrt{n}} \sum_{i=1}^{n}\left\{(f_1(X_i), \ldots ,
  f_k(X_i))-E(f_1(X_0), \ldots , f_k(X_0))\right\}
  \cla N(0,\Sigma_{f_1,\ldots,f_k}),
\end{equation*}
where $N(0,\Sigma_{f_1,\ldots,f_k})$ denotes a multivariate normal law with 
mean zero and covariance matrix 
\[ 
 \Sigma_{f_1,\ldots, f_k} = (\sigma_{f_i},_{f_j})_{1\leq i, j\leq k}
\]
where for any Lipschitz functions $ f, g $ we define
\begin{eqnarray*}
 \sigma_{f,g}= \Cov (f(X_0),g (X_0))
 &+& \sum^\infty_{k=1} \Cov (f(X_0),g(X_k))\\
 &+& \sum^\infty_{k=1} \Cov (f(X_k),g(X_0)).
\end{eqnarray*}

This result proves the proposition. \hfill $\Box$

\begin{proposition} 
For any $\varepsilon,\eta > 0$ there exists a partition $0=t_0^\prime<\ldots<t_m^\prime=1$ 
such that
\[
  \limsup_{n\rightarrow\infty} 
  P \left(\sup\limits_{0\leq t\leq 1}\left| U_n(t) - U_n^{(m)}(t)\right|>
 \varepsilon\right) \leq \eta.
\]
\label{prop:ep-appr}
\end{proposition}

{\em Proof.}

By a variant of the well known chaining technique we will control
\[
  P\left( \sup\limits_{0 \leq t \leq 1} \left| U_n (t) - U_n^{(m)}(t)\right| 
  \geq \varepsilon\right),
\]
and then show that this probability can be made arbitrarily small by choosing 
a partition $0=t_0^\prime < \ldots < t_m^\prime=1$ that is fine enough.
From here on we assume that the partition $0=t_0^\prime < \ldots < t_m^\prime=1$ is regularly distributed.
Let $h=\frac{1}{m}=t_j^\prime-t_{j-1}^\prime$, for $j=1,\dots,m$.

On the interval $[t_{j-1}^\prime, t_j^\prime]$ we introduce a sequence of refining partitions
\[
 t_{j-1}^\prime = s_0^{\prime(k)} < s_1^{\prime(k)} < \ldots < s^{\prime(k)}_{2^k} = t_j^\prime
\]
by
\[
 s_l^{\prime(k)} = t^\prime_{j-1} + l \cdot \frac{h}{2^k}\quad,
 \quad 0 \leq l \leq 2^k.
\]
Let us define
$$
s_l^{(k)}=F^{-1}(s_l^{\prime(k)})\quad,\quad 0 \leq l \leq 2^k.
$$
We now have partitions of $[t_{j-1},t_j]$,
\[
 t_{j-1} = s_0^{(k)} < s_1^{(k)} < \ldots < s^{(k)}_{2^k} = t_j.
\]
For convenience, we also consider the points
$$
s_{-1}^{(k)}=F^{-1}\left(t^\prime_{j-1} - \frac{h}{2^k}\right)
$$
and the points
$$
s_{2^k+1}^{(k)}=F^{-1}\left(t^\prime_{j-1} + (2^k+1)\frac{h}{2^k}\right).
$$

For any $t\in [t_{j-1}, t_j)$ and $k\geq 0$ we define the index
\[
 l(k,t) = \max \left\{l: s_l^{(k)} \leq t \right\}.
\]
In this way we obtain a chain
\[
  t_{j-1} = s_{l(0,t)}^{(0)} \leq s_{l(1,t)}^{(1)} \leq \ldots
  \leq s_{l(k,t)}^{(k)} \leq t \leq s_{l(k,t)+1}^{(k)},
\]
linking the left endpoint $t_{j-1}$ to $t$. Note that for 
$t\in [t_{j-1}, t_j)$ we have by definition
$U_n^{(m)}(t) = U_n^{(m)}(t_{j-1})$. We define the functions $\psi^{(k)}_l$,
$k \geq 0$, $0 \leq l \leq 2^k$, by
\[
 \psi^{(k)}_l(x) = \varphi\left(\frac{x}{s_{l}^{(k)}-s_{l-1}^{(k)}}\right),
\]
where $\varphi$ is defined as in (\ref{eq:phi}). Note that 
$\psi^{(0)}_{l(0,t)}(x-s^{(0)}_{l(0,t)})=\varphi_j(x)$. To be consistent, in the case $j=1$, we have to fix
$\psi_0^{(k)}\equiv 0$, for all $k\ge 0$.
We build a chain bridging the gap between
\[
 F_n(t) = \frac{1}{n} \sum\limits^n_{i=1} 1_{(-\infty, t]}(X_i)
\]
and
\[
 F_n^{(m)}(t) = \frac{1}{n} \sum\limits^n_{i=1} \varphi_j (X_i)
\]
by the functions
\begin{eqnarray*}
 \varphi_j(x)
 &=&\psi^{(0)}_{l(0,t)}( x-s^{(0)}_{l(0,t)}) \\
 &\leq & \psi^{(1)}_{l(1,t)} (x-s^{(1)}_{l(1,t)}) \\
 &\leq& \ldots \\
 &\leq& \psi^{(K)}_{l(K,t)}(x-s^{(K)}_{l(K,t)}) \\
 &\leq&  1_{(-\infty, t]}(x) \\
 &\leq & \psi^{(K)}_{l(K,t)+2}(x-s^{(K)}_{l(K,t)+2}),
\end{eqnarray*}
where $K$ is some integer to be chosen later. In this way we get
\begin{eqnarray}
 F_n(t) - F_n^{(m)}(t) 
  &=& \sum^K_{k=1} \frac{1}{n} 
   \sum^n_{i=1} \left(\psi^{(k)}_{l(k,t)}(X_i - s^{(k)}_{l(k,t)})
       -\psi^{(k-1)}_{l(k-1,t)}(X_i-s^{(k-1)}_{l(k-1,t)})\right)\nonumber \\
  && + \frac{1}{n}\sum^n_{i=1}\left( 1_{(-\infty, t]}(X_i)
    - \psi^{(K)}_{l(K,t)}(X_i-s^{(K)}_{l(K,t)})\right).
\label{eq:fn-fnm}
\end{eqnarray}
Observe that by definition of $s^{(k)}_{l(k,t)}$ and of $\psi^{(K)}$,
\begin{eqnarray*}
 0 &\leq & 1_{(-\infty,t]} (X_i) - \psi^{(K)}_{l(K,t)} ( X_i - s^{(K)}_{l(K,t)}
 ) \\
  & \leq & \psi^{(K)}_{l(K,t)+2} (X_i - s^{(K)}_{l(K,t)+2}) - \psi^{(K)}_{l(K,t)}
  (X_i - s^{(K)}_{l(K,t)}).
\end{eqnarray*}

From (\ref{eq:fn-fnm}) we get by centering and normalization
\begin{eqnarray*}
 U_n(t) - U_n^{(m)}(t)
 &=& \sum^K_{k=1} 
 \frac{1}{\sqrt{n}}\sum^n_{i=1}
  \left\{ 
 \left( \psi^{(k)}_{l(k,t)}( X_i-s^{(k)}_{l(k,t)})-E \psi^{(k)}_{l(k,t)}(X_i-s^{(k)}_{l(k,t)} )
 \right) \right.  \\
 && -\left. \left(\psi^{(k-1)}_{l(k-1,t)}(X_i - s^{(k-1)}_{l(k-1,t)})
  -E\psi^{(k-1)}_{l(k-1,t)}(X_i - s^{(k-1)}_{l(k-1,t)})\right)\right\} \\
 && + \frac{1}{\sqrt{n}}  \sum^{n}_{i=1}
  \left\{\left(1_{(-\infty,t]}(X_i)-F(t)\right)\right.\\
 && - \left. \left(\psi^{(K)}_{l(K,t)}(X_i - s^{(K)}_{l(K,t)}) 
 -E \psi^{(K)}_{l(K,t)}(X_i - s^{(K)}_{l(K,t)})\right)\right\}.
\end{eqnarray*}
For the last term on the r.h.s. we have the following upper and lower bounds,
\begin{eqnarray*}
 &&
 \frac{1}{\sqrt{n}} \sum^n_{i=1} \left\{ 
 \left(1_{(-\infty, t]} (X_i)-F(t)\right) - 
 \left(\psi^{(K)}_{l(K,t)} (X_i - s^{(K)}_{l (K,t)}) - E \psi^{(K)}_{l(K,t)}
 (X_i - s_{l (K,t)}^{(K)})\right)\right\} \\
&& \le  \frac{1}{\sqrt{n}} \sum^n_{i=1} \left\{  
 \left( \psi^{(K)}_{l(K,t)+2} (X_i-s^{(K)}_{l
 (K,t)+2}) - E \psi^{(K)}_{l(K,t)+2} (X_i- s^{(K)}_{l (K,t)+2} ) \right)\right. \\
&& \quad \quad -\left. \left( \psi^{(K)}_{l(K,t)} (X_i - s^{(K)}_{l (K,t)} ) 
 - E \psi^{(K)}_{l(K,t)}(X_i - s^{(K)}_{l (K,t)} ) \right) \right\}\\
&&\quad +\sqrt{n} \left( E\psi^{(K)}_{l(K,t)+2} (X_i - s^{(K)}_{l(K,t)+2})-F(t)\right)
\end{eqnarray*}
and
\begin{eqnarray*}
 && \frac{1}{\sqrt{n}}\sum^{n}_{i=1}  \left\{ 
 \left( 1_{(-\infty,t]}(X_i) - F(t)\right)-
 \left(\psi^{(K)}_{l(K,t)}(X_i -s^{(K)}_{l(K,t)}) - E \varphi^{(K)}_{l(K,t)} 
 (X_i - s^{(K)}_{l(K,t)})\right)\right\} \\
 && \geq - \sqrt{n}\left(F(t) - E \psi^{(K)}_{l(K,t)} 
  (X_i - s^{(K)}_{l(K,t)})\right).
\end{eqnarray*}
Now choose $K = 4 + \left\lfloor \log\left( \frac{\sqrt{n}h}{\varepsilon}\right)\log^{-1}(2)\right\rfloor$
 and note that 
\begin{equation*}
\frac{\varepsilon}{2^4}\leq\sqrt {n} \frac{h}{2^K} \leq \frac{\varepsilon}{2^3}
\end{equation*}
and thus
\begin{eqnarray*}
&&\sqrt{n}\left| E \psi^{(K)}_{l(K,t)+2}(X_i - s^{(K)}_{l(K,t)+2} )- 
E \psi^{(K)}_{l(K,t)}(X_i -  s^{(K)}_{l(K,t)})\right|\\
&&\leq\sqrt{n}\left|F(s^{(K)}_{l(K,t)+2})-F(s^{(K)}_{l(K,t)-1})\right|\\
&&\leq\frac{\varepsilon}{2}.
\end{eqnarray*}
Thus we get for all $t \in [t_{j-1}, t_j]$, 
\begin{eqnarray*}
 && \hspace{-15mm}\left| U_n(t)-U_n^{(m)} (t)\right| \\
 &\leq& \sum^K_{k=1} \frac{1}{\sqrt{n}} 
  \left| \sum^n_{i=1}
  \left\{  \left( \psi^{(k)}_{l(k,t)} (X_i-s^{(k)}_{l(k,t)} ) - 
  E \psi^{(k)}_{l(k,t)} (X_i - s^{(k)}_{l (k,t)}) \right) \right.\right. \\
&& \left.\left.\quad - \left( \psi^{(k-1)}_{l(k-1,t)}  (X_i-s^{(k-1)}_{l (k-1,t)}) - 
  E \psi^{(k-1)}_{l(k-1,t)}  (X_i -
  s^{(k-1)}_{l (k-1,t)})\right)  \right\} \right| 
  \\
&& + \frac{1}{\sqrt{n}} \left|  \sum\limits^n_{i=1} \left\{  
 \left( \psi^{(K)}_{l(K,t)+2}  (X_i -
 s^{(K)}_{l (K,t)+2}) - E \psi^{(K)}_{l(K,t)+2}  (X_i - s^{(K)}_{l (K,t)+2})\right)\right.\right. \\
&&\left.\left. \quad -\left( \psi^{(K)}_{l(K,t)}  (X_i - s^{(K)}_{l (K,t)} ) - E \psi^{(K)}_{l(K,t)} 
  (X_i - s^{(K)}_{l (K,t)}) \right)   \right\} 
   \right| \\
&& \quad + \frac{\varepsilon}{2}.
\end{eqnarray*}
Note that by definition of $l(k,t)$ and of $s_{l}^{(k)}$, we have
$ s^{(k-1)}_{l (k-1,t)} \in \{ s^{(k)}_{l (k,t)}, s^{(k)}_{l (k,t)-1}\}$
and thus
$$
l(k-1,t)=\left\lfloor \frac{l(k,t)}{2} \right\rfloor.
$$
Therefore
\begin{eqnarray*}
&& \hspace{-15mm} \sup_{t_{j-1} \le t \le t_j}  \left| 
    U_n(t) - U^{(m)}_n (t) \right| \\
&\le&  \sum^K_{k=1} \frac{1}{\sqrt{n}} \max_{0 \le l \le 2^k-1}
   \left|  \sum^n_{i=1} 
  \left( ( \psi^{(k)}_l (X_i - s_{l}^{(k)} ) - E \psi^{(k)}_l
  (X_i - s_{l}^{(k)}) )\right. \right. \\
&& \qquad \quad \left. \left. - ( \psi^{(k-1)}_{\lfloor\frac{l}{2}\rfloor} (X_i - s_{\lfloor\frac{l}{2}\rfloor}^{(k-1)}) - E \psi^{(k-1)}_{\lfloor\frac{l}{2}\rfloor}
  (X_i - s_{\lfloor\frac{l}{2}\rfloor}^{(k-1)} )) \right)  \right|\\
&& \quad + \frac{1}{\sqrt{n}} \max_{0\le l \le 2^K-1}   \left|  
  \sum^n_{i=1} \left( (\psi^{(K)}_{l+2} (X_i - s_{l + 2}^{(K)}) - E \psi^{(K)}_{l+2} 
  (X_i - s^{(K)}_{l+2})) \right. \right. \\
&& \left.\left. \qquad \quad - ( \psi^{(K)}_l (X_i - s_{l}^{(K)}) - E \psi^{(K)}_l
  (X_i - s_l^{(K)} )) \right) \right| \\
&& \quad + \frac{\varepsilon}{2}.
\end{eqnarray*}

\medskip

Now take $\varepsilon_k  : = \frac{\varepsilon}{4 k (k+1)} $ and note that $\sum^K_{k=1}\varepsilon_k \leq 
\frac{\varepsilon}{4}.$

Then we obtain
\begin{eqnarray*}
&&\hspace{-15mm} P \left(\sup_{t_{j-1}\le t \le t_j} 
  \left|U_n (t) - U_n^{(m)} (t) \right| \ge\varepsilon\right)\\
& \le &\sum\limits^K_{k=1} \sum\limits^{2^k-1}_{l =0} 
P \left( \frac{1}{\sqrt{n}}\right. \left| \sum\limits^n_{i=1} \right. \left\{
\left( \psi^{(k)}_l (X_i - s_{l}^{(k)}) - E \psi^{(k)}_l 
(X_i - s_l^{(k)}) \right) \right. \\
&& \quad \left.\left.\left.-\left( \psi^{(k-1)}_{\lfloor\frac{l}{2}\rfloor} (X_i-s_{\lfloor\frac{l}{2}\rfloor}^{(k-1)}) - E \psi^{(k-1)}_{\lfloor\frac{l}{2}\rfloor} 
 (X_i-s_{\lfloor\frac{l}{2}\rfloor}^{(k-1)})
\right) \right\} \right| \ge \varepsilon_k \right) \\
&& + \sum^{2^K-1}_{l =0} P \left( \frac{1}{\sqrt{n}}\right.  \left| 
\sum^n_{i=1}\left\{ \left( \psi^{(K)}_{l+2} (X_i-s_{l +2}^{(K)}) - 
E \psi^{(K)}_{l+2} (X_i -s_{l + 2}^{(K)} ) \right)\right.\right. \\
&& \quad \left.\left.\left.-\left(\psi^{(K)}_l(X_i-s_{l}^{(K)}) - E \psi^{(K)}_l (X_i-s_{l}^{(K)})
\right)  \right\} \right|  \ge \frac{\varepsilon}{4}\right).
\end{eqnarray*}
At this point we use Markov's inequality together with the 4-th moment 
bound (\ref{eq:4th-moment}).
\begin{eqnarray*}
&&\hspace{-15mm} P \left(\sup_{t_{j-1}\le t \le t_j} 
  \left|U_n (t) - U_n^{(m)} (t) \right| \ge\varepsilon\right)\\
& \le & C\sum\limits^K_{k=1} \sum\limits^{2^k-1}_{l =0} 
\left\{\frac{1}{n\varepsilon_k^4}\left\| \psi^{(k)}_l(X_0 - s_l^{(k)}) - \psi^{(k-1)}_{\lfloor\frac{l}{2}\rfloor} 
 (X_0 - s_{\lfloor\frac{l}{2}\rfloor}^{(k-1)}) \right\|_1 \right.\\
&&\hspace{60mm}.\log^\alpha\left(1+\left\|\psi^{(k)}_l - \psi^{(k-1)}_{\lfloor\frac{l}{2}\rfloor} \right\|\right)\\
&&+ \frac{1}{\varepsilon_k^4}\left\| \psi^{(k)}_l(X_0 - s_l^{(k)}) - \psi^{(k-1)}_{\lfloor\frac{l}{2}\rfloor} 
 (X_0 - s_{\lfloor\frac{l}{2}\rfloor}^{(k-1)}) \right\|_1^2\\
&&\hspace{60mm}\left. .\log^\beta\left(1+\left\|\psi^{(k)}_l - \psi^{(k-1)}_{\lfloor\frac{l}{2}\rfloor} \right\|\right)\right\}\\
&&+C\sum\limits^{2^k-1}_{l =0} 
\left\{\frac{4^4}{n\varepsilon^4}\left\| \psi^{(K)}_{l+2}(X_0 - s_{l+2}^{(K)}) - \psi^{(K)}_{l} 
 (X_0 - s_{l}^{(K)}) \right\|_1\right. \\
&&\hspace{60mm}.\log^\alpha\left(1+\left\|\psi^{(K)}_{l+2} - \psi^{(K)}_{l} \right\|\right)\\
&&+ \frac{4^4}{\varepsilon^4}\left\| \psi^{(K)}_{l+2}(X_0 - s_{l+2}^{(K)}) - \psi^{(K)}_{l} 
 (X_0 - s_{l}^{(K)}) \right\|_1^2\\
&&\hspace{60mm}\left. .\log^\beta\left(1+\left\|\psi^{(K)}_{l+2} - \psi^{(K)}_{l} \right\|\right)\right\}.
\end{eqnarray*}
Note that
\begin{eqnarray*}
 \left\| \psi^{(k)}_l(X_0 - s_l^{(k)}) - \psi^{(k-1)}_{\lfloor\frac{l}{2}\rfloor} 
 (X_0 - s_{\lfloor\frac{l}{2}\rfloor}^{(k-1)}) \right\|_1 
&\leq & \left| F(s_l^{(k)}) - F(s_{{\lfloor\frac{l}{2}\rfloor}-1}^{(k-1)}) \right|\\
&\leq & \left| F(s_l^{(k)}) - F(s_{l-3}^{(k)}) \right|\\
& =&  \frac{3h}{2^k}
\end{eqnarray*}
and
\begin{eqnarray*}
 \left\| \psi^{(K)}_{l+2}(X_0 - s_{l+2}^{(K)}) - \psi^{(K)}_{l} 
 (X_0 - s_{l}^{(K)}) \right\|_1&\le &
\left| F(s_{l+2}^{(K)}) - F(s_{l-1}^{(K)}) \right|\\
& =&  \frac{3h}{2^K}.
\end{eqnarray*}

\medskip

If (\ref{eq-modulus}) is satisfied, 
\begin{eqnarray*}
 \left\| \psi^{(k)}_l \right\|
&\leq& 1+\left[\inf\left\{s>0 : \forall t, F(t+s)-F(t)\geq \frac{h}{2^k}\right\}\right]^{-1}\\
&\leq& 1+\left[\inf\left\{s>0 : D|\log(s)|^{-\gamma}\geq \frac{h}{2^k}\right\}\right]^{-1}\\
&=& 1+  \exp\left(\left(\frac{D2^k}{h}\right)^\frac{1}{\gamma}\right).
\end{eqnarray*}

Thus we have
\begin{eqnarray*}
  && \hspace*{-20mm} P \left(\sup\limits_{ t_{j-1} \leq t \leq t_j}
 \left| U_n (t) - U_n (t_j)\right|\geq \varepsilon\right) \\
 &\leq& 4^4 C \sum\limits^K_{k=1} 2^k \frac{(k(k+1))^4}{\varepsilon^4} 
 \frac{1}{n} \frac{3h}{2^k}\log^\alpha \left(2+\exp\left(\left(\frac{D2^k}{h}\right)^\frac{1}{\gamma}\right)
  \right)\\
 && + 4^4 C\sum\limits^K_{k=1} 2^k \frac{(k(k+1))^4}{\varepsilon^4}
 \frac{(3h)^2} {2^{2k}} \log^\beta \left(2+\exp\left(\left(\frac{D2^k}{h}\right)^{\frac{1}{\gamma}}\right)\right)\\
 && + 4^4 C2^K \frac{1}{\varepsilon^4}\frac{1}{n}\frac{3h}{2^K}
  \log^\alpha\left(2+ \exp\left(\left(\frac{D2^k}{h}\right)^{\frac{1}{\gamma}}\right)\right)\\
 && + 4^4 C2^K \frac{1}{\varepsilon^4} \frac{(3h)^2}{2^{2K}} 
 \log^\beta\left(2+ \exp\left(\left(\frac{D2^k}{h}\right)^{\frac{1}{\gamma}}\right)\right)\\
 &\leq& \frac{1}{n} \frac{C^\prime}{\varepsilon^4}
  \sum^K_{k=1} k^8 h
  \left(\frac{D2^k}{h}\right)^{\frac{\alpha}{\gamma}}
  + \frac{C^\prime}{\varepsilon^4}
 \sum^K_{k=1}\frac{k^8}{2^k}h^2 
 \left(\frac{D2^k}{h}\right)^{\frac{\beta}{\gamma}}\\ 
 &\leq& D^{\frac{\alpha}{\gamma}}\frac{1}{n} \frac{C^\prime}{\varepsilon^4}h
 \left(\frac{2^K}{h}\right)^{\frac{\alpha}{\gamma}}
 \sum^K_{k=1} k^8
   + D^{\frac{\beta}{\gamma}}\frac{C^\prime}{\varepsilon^4}h^{2-\frac{\beta}{\gamma}}
 \sum^\infty_{k=1}k^8 2^{k(\frac{\beta}{\gamma}-1)}\\
&\leq&\frac{h}{n}\frac{C''}{\varepsilon^4}\left(\frac{\sqrt{n}}{\varepsilon}\right)^{\frac{\alpha}{\gamma}}K^9
 + \frac{C''}{\varepsilon^4}h^{2-\frac{\beta}{\gamma}}
\end{eqnarray*}
where $C'$ and $C''$ are some constants and we have used convergence of the series
$\sum^\infty_{k=1}k^8 2^{k(\frac{\beta}{\gamma}-1)}$.

\medskip

Finally, using $mh=1$,
\begin{eqnarray*}
&& \hspace*{-20mm} P \left( \sup_{0\leq t \leq 1} 
 \left| U_n (t) - U_n^{(m)}(t) \right| 
\geq \varepsilon \right) \\
 &\leq& \sum\limits^m_{j=1} P \left( \sup_{t_{j-1} \leq t \leq t_j} 
 \left| U_n(t) - U_n^{(m)} (t) \right| \geq \varepsilon  \right)\\
 &\leq& m h n^{\frac{\alpha}{2\gamma}-1}\frac{C''}{\varepsilon^{4+\frac{\alpha}{\gamma}}}K^9
 + m\frac{C''}{\varepsilon^4}h^{2-\frac{\beta}{\gamma}}\\
&\leq& n^{\frac{\alpha}{2\gamma}-1}\frac{C''}{\varepsilon^{4+\frac{\alpha}{\gamma}}}
\left(4+\log\frac{\sqrt{n}h}{\varepsilon}\right)^9
+ \frac{C''}{\varepsilon^4}h^{1-\frac{\beta}{\gamma}}
\end{eqnarray*}
Now, the first of the two final summands converges to zero as 
$n \rightarrow \infty$. The second can be made arbitrarily  
small by choosing a partition that is fine enough (i.e. $h$ small). \hfill $\Box$

\bigskip

We used a different technique than the usual finite dimensional convergence plus tightness. Of course, since the weak convergence implies the finite dimensional convergence and the tightness, these two properties are satified.
Nevertheless, we can also deduce a tightness criterion implying that, almost surely, the limit process has continuous sample paths (see Billingsley \cite{Bil68}, Theorem 15.5).

\medskip

\begin{proposition}
For all $\varepsilon,\, \eta>0$, there exist $\delta>0$ and $N\ge 0$ such that for all $n\ge N$,
$$
P\left(\sup_{|t-s|<\delta}|U_n(t)-U_n(s)|\ge \varepsilon\right)\le \eta.
$$
In particular,  $P(W\in C(\R))=1$.
\end{proposition}

{\em Proof.}

Let $\varepsilon>0$ and $\eta>0$.
Let $m$ be an integer such that 
\begin{equation}\label{m}
\frac{C}{\varepsilon^4}\frac{D^{\frac{\beta}{\gamma}}}{m^{1+\frac{\beta}{\gamma}}}<\frac{\eta}{4}
\end{equation}
and consider the regular partition of $[0,1]$ with mesh $\frac{1}{m}$.

By Proposition \ref{prop:ep-appr}, there exists $N\ge 0$ such that for all
$n\ge N$,
$$
P \left(\sup\limits_{0\leq t\leq 1}\left| U_n(t) - U_n^{(m)}(t)\right|\ge
 \frac{\varepsilon}{3}\right) \leq \frac{\eta}{4}.
$$ 
Let $\delta>0$ such that $\delta<\frac{1}{m}$. Then, for all $n\ge N$,
\begin{eqnarray*}
&&\hspace{-20pt} P\left(\sup_{|t-s|<\delta}|U_n(t)-U_n(s)|\ge \varepsilon\right)\\
&&\le2P\left(\sup_{0\le t\le 1}|U_n(t)-U_n^m(t)|\ge \frac{\varepsilon}{3}\right)
+P\left(\sup_{|t-s|<\delta}|U_n^m(t)-U_n^m(s)|\ge \frac{\varepsilon}{3}\right)\\
&&\le \frac{\eta}{2}+P\left(\sup_{|t-s|<\delta}|U_n^m(t)-U_n^m(s)|\ge \frac{\varepsilon}{3}\right).
\end{eqnarray*}
We recall, as $t_j=F^{-1}(t_j')=F^{-1}(\frac{j}{m})$, that
\begin{eqnarray*}
 \|\varphi_j(X_0)-\varphi_{j+1}(X_0)\|_1&\le&P\left(t_{j-2}\le X_0\le t_j\right)\le \frac{2}{m},\\
\|\varphi_j\| & \le & 1+\exp\left(\left(\frac{D}{m}\right)^{\frac{1}{\gamma}}\right).
\end{eqnarray*}
Thus, by the 4-th moment bound (\ref{eq:4th-moment}), 
\begin{eqnarray*}
P\left(\sup_{|t-s|<\delta}|U_n^m(t)-U_n^m(s)|\ge \frac{\varepsilon}{3}\right)
\le \frac{C}{n\varepsilon^4}\left(\frac{D}{m}\right)^{\frac{\alpha}{\gamma}}
+\frac{C}{\varepsilon^4}\frac{D^{\frac{\beta}{\gamma}}}{m^{1+\frac{\beta}{\gamma}}}.
\end{eqnarray*}
Now there exists $N'\ge N$ such that 
$$
\frac{C}{n\varepsilon^4}\left(\frac{D}{m}\right)^{\frac{\alpha}{\gamma}}\le \frac{\eta}{4}.
$$
Finally, by (\ref{m}),
$$
P\left(\sup_{|t-s|<\delta}|U_n(t)-U_n(s)|\ge \varepsilon\right)\le \eta.
$$
\hfill $\Box$

\

\subsection{The unbounded case}

Let $(X_i)_{i\ge 0}$ be an $\R$-valued stationary ergodic random process such that
(\ref{eq:lip-clt}), (\ref{eq:4th-moment}) and (\ref{eq-modulus}) hold. We will
show that it can be reduced to the case of bounded variables.

\medskip

For all $x<y\in\R$, we say that the closed interval $[x,y]$ is a 'bad' interval (for $F$) if
$$
F(y)-F(x)\ge y-x.
$$
We say that $[x,y]$ is a maximal 'bad' interval (for $F$) if for all 'bad' intervals
$[a,b]$, we have $[a,b]\subset[x,y]$ or $[a,b]\cap[x,y]=\emptyset$.

We denote by $I^{max}$ the set of all maximal 'bad' intervals.

\begin{lemma} 
$ $

\begin{itemize}
\item[(i)] The Lebesgue measure of 
$$I:=\bigcup_{[x,y]\in I^{max}}[x,y]$$
is smaller than $1$.
\item[(ii)] For all $[x,y]\in I^{max}$, we have
$$
F(y)-F(x)=y-x.
$$ 
 \end{itemize}
\end{lemma}

{\em Proof.}

Because $F$ is non-decreasing and takes values in $[0,1]$, the first assertion is clear.

\medskip

If for $x<y$, $F(y)-F(x)>y-x$, then there exists $\varepsilon>0$ such that
$$
F(y)-F(x)>y-x+\varepsilon.
$$
Thus, for all
$z>y$ such that $z-y\le\varepsilon$, by monotonicity of $F$, we have
\begin{eqnarray*}
 F(z)-F(x)&\ge& F(y)-F(x)\\
&>&y-x+\varepsilon\\
&\ge&z-x
\end{eqnarray*}
and then $[x,y]$ is not maximal.\hfill $\Box$

\medskip

We define the function $g$ from $\R$ to $]0,1[$ by 
$$
\mbox{ for all }[x,y]\in I^{max}, \mbox{ for all }t\in [x,y],\, g(t):=F(x)+t-x
$$
and
$$
\mbox{ for all }t\notin I,\, g(t):=F(t).
$$
Then $g$ is a 1-Lipschitz function.

\medskip

We define the $[0,1]$-valued stationary ergodic random process
$(Y_i)_{i\ge 0}$ by
$$
Y_i=g(X_i),\;i\ge 0.
$$
Since $g$ is Lipschitz, $(Y_i)_{i\ge 0}$ satisfies (\ref{eq:lip-clt}) and (\ref{eq:4th-moment}).

We also have
\begin{equation*}
 G(t):=P(Y_0\le t)=F\circ g^{-1}(t)
\end{equation*}
where 
$$
g^{-1}(t)=\sup\{s\in\R\,:\,F(s)\le t\}.
$$
Clearly, $G$ is the identity on $g(\R\setminus I)$. Further, for all $[x,y]\in I^{max}$, the graph
of $G$ on $g([x,y])$ is the graph of $F$ on $[x,y]$ and the Lebesgue measure of $g([x,y])$ is equal to the Lebesgue measure of $[x,y]$.
Then 
\begin{eqnarray*}
 \omega_G(\delta)&\le &\max\{\omega_F(\delta),\delta\}
\end{eqnarray*}
and (\ref{eq-modulus}) holds.

\medskip

We define the associated distribution functions and empirical processes
\begin{eqnarray*}
 F_n(t) &:= & \frac{1}{n} \sum_{i=1}^n 1_{(-\infty,t]} (X_i),\;  t\in\R,
 \\
 U_n(t) &:=& \sqrt{n} (F_n(t)-F(t)),\; t \in\R,
\\
 G_n(t) &:= & \frac{1}{n} \sum_{i=1}^n 1_{[0,t]} (Y_i),\; 0\leq t\leq 1,
 \\
 V_n(t) &:=& \sqrt{n} (G_n(t)-G(t)),\; 0\leq t \leq 1.
\end{eqnarray*}
We have 
$$
U_n(t)=V_n(g(t)),\; t \in\R.
$$
By the theorem for bounded variables (section \ref{boundedcase}), 
\begin{equation*}
 (V_n(t))_{0\leq t \leq 1} \stackrel{\mathcal{D}}
 {\longrightarrow} (V(t))_{0\leq t \leq1},
\end{equation*}
where $V(t)$ is a mean-zero Gaussian process such that $P(V\in C[0,1])=1$.

Applying Theorem 5.1 of Billingsley \cite{Bil68} with
\begin{eqnarray*}
 h:D[0,1]&\longrightarrow& D(\R)\\
x&\mapsto& x\circ g,
\end{eqnarray*}
we get the weak convergence of $(U_n(t))_{ t \in \R}$ to a Gaussian process 
$$
(W(t))_{t\in\R}=(V\circ g(t))_{t\in\R}
$$
such that $P(W\in C(\R))=1$.

\section{Proof of Theorem \ref{pbil}}

\begin{lemma}
Let $(X,d)$ be a complete metric space and let $x_n, x_n^{(m)}, x^{m} 
\in X$, $n \geq 1, m \geq 1$ be given with the properties
\begin{eqnarray}
 \lim_{n\rightarrow \infty} d (x_n^{(m)}, x^{(m)}) &=& 0 \qquad \forall m 
 \label{eq:l1} \\
 \lim_{m\rightarrow \infty} \limsup_{n\rightarrow \infty} d (x_n, x_n^{(m)}) 
 &=& 0 .\label{eq:l2}
\end{eqnarray}
Then $x:= \lim_{m\rightarrow \infty} x^{(m)}$ exists and 
\[
 \lim_{n\rightarrow \infty} d (x_n, x) = 0.
\]
\end{lemma}
{\em Proof.}
We will first show that $x^{(m)}$ is a Cauchy sequence. 
Given $ \epsilon > 0,$ choose $M$ so big that $ \forall m \geq M $
\[
 \limsup_{n\rightarrow \infty}d(x_n, x_n^{(m)})< \frac{\varepsilon}{4}. 
\]
Now take $m_1,m_2 \geq M$. For all $n$ sufficiently large, we have then
\begin{eqnarray*}
 d (x_n^{(m_1)}, x^{(m_1)}) &<& \frac{\varepsilon}{4} \\
 d (x_n^{(m_2)}, x^{(m_2)}) &<& \frac{\varepsilon}{4} \\
 d (x_n, x_n^{(m_1)}) &<& \frac{\varepsilon}{4} \\
 d (x_n, x_n^{(m_2)}) &<& \frac{\varepsilon}{4},
\end{eqnarray*}
and hence, by the triangle inequality $d (x^{(m_1)}, x^{(m_2)}) < 
\varepsilon.$ Thus $(x^{(m)})_{m \geq 1}$
is a Cauchy sequence and hence $x: = \lim_{m\rightarrow \infty} x^{(m)} $ 
exists.
\\
It remains to show that $ \lim_{n\rightarrow \infty} x_n = x. $
Given $\varepsilon > 0, $ choose $m_0 $ so that
\[
  \limsup_{n\rightarrow\infty} d (x_n, x_n^{(m_0)}) 
 <\frac{\varepsilon}{4} 
\]
and $d (x^{(m_0)}, x)< \frac{\varepsilon}{4}$. 
Then choose $N $ such that for all $n \geq N $
\begin{eqnarray*}
 d (x_n, x_n^{(m_0)}) &<& \frac{\varepsilon}{4}\\
 d (x_n^{(m_0)}, x^{(m_0)}) &<& \frac{\varepsilon}{4}.
\end{eqnarray*}
Using the triangle inequality, we get
\[d(x_n, x) <\varepsilon 
\]
for all $n \geq N. $ \hfill $\Box$ 
\vspace*{3mm}

{\em Proof of Theorem \ref{pbil}.}
Let $\mu_n, \mu_n^{(m)}$ and $\mu^{(m)}$ denote the distributions of the 
random variables $X_n, X_n^{(m)}$ and $X^{(m)}$ respectively. These are 
elements of $M_1 (S),$ the space of probability measures on $S$. 
We consider the Prohorov metric $d$ on $ M_1 (S)$, defined by
\[ 
 d (\mu, \upsilon) = \inf \left\{ \varepsilon > 0 : \mu (A) 
 \leq \upsilon(A^\varepsilon) + \varepsilon\quad \forall A \subset 
 S\; \text{measurable} \right\}. 
\]
Note that $ (M_1(S),d)$ is a complete metric space. 
If $ Y, Z$ are two S-valued random variables with distributions 
$P_Y, P_Z,$ satisfying
\[  
  P( \rho ( Y,Z) \geq \varepsilon ) \leq \varepsilon, 
\]
then $d(P_Y, P_Z) \leq \varepsilon.$
Moreover $d$ metrizes the topology of weak convergence, i. e. 
$\mu_n \rightarrow \mu$ if and only if $d (\mu_n, \mu) \rightarrow 0.$
We now apply Lemma 3.1 to $ \mu_n, \mu_n^{(m)}, \mu^{(m)}.$
Note that (\ref{eq:l1}) is a direct consequence of (\ref{eq:p1}). 
Given $\varepsilon > 0,$ by (\ref{eq:p2}) we can find $m_0$ 
such that for all $m \geq m_0,$
\[ 
 \limsup_{n\rightarrow \infty} P (\rho (X_n, X_n^{(m)}) 
 \geq \varepsilon ) < \varepsilon.
\]
Fix such an $m$; then we can find $n_0$ such that $\forall n \geq n_0$
\[ 
 P ( \rho (X_n, X_n^{(m)}) \geq \varepsilon ) \leq \varepsilon
\]
and thus $d (\mu_n, \mu_n^{(m)}) \leq \varepsilon.$
Hence 
\[
 \limsup_{n \rightarrow \infty} d (\mu_n, \mu_n^{(m)}) \leq \varepsilon 
\]
for all $m \geq m_0,$ showing that (\ref{eq:l2}) holds. 
Thus by Lemma 3.1, there exists a probability distribution $\mu $ on $S$ 
such that
\begin{eqnarray*}
 \lim_{n\rightarrow \infty} d (\mu^{(m)}, \mu) &=& 0 \\
 \lim_{m\rightarrow \infty} d (\mu_n , \mu) &=& 0 .
\end{eqnarray*}
Finally, let $X$ be an S-valued random variable with distribution 
$\mu$. Then $X^{(m)}\stackrel {\mathcal{D}}\rightarrow X$ as 
$m \rightarrow \infty$ and $X_n \rightarrow X$ as $n \rightarrow \infty.$
\hfill $\Box$

\section{Examples}\label{examp}

According to Durieu \cite{Dur08b} the 4-th moment bound (\ref{eq:4th-moment}) holds for Markov chains and dynamical systems
under some assumptions on the Markov transition operator or the Perron-Frobenius operator.

\medskip

Let $(E,d)$ be a separable metric space and $(X_k)_{k\ge0}$ be an $E$-valued Markov chain with transition operator $Q$ and invariant measure $\nu$.
Denote by $\LL$ the space of all bounded Lipschitz continuous functions from $E$ to $\R$ equipped with the norm defined in
(\ref{eq:4th-moment}). We say that the Markov chain $(X_k)_{k\ge0}$ is $\LL$-geometrically ergodic if there exist $C>0$ and $0<\theta<1$ such that for all $f\in\LL$,
\begin{equation}\label{geom}
\|Q^kf-\Pi f\|\le C\theta^k\|f\|,
\end{equation}
where $\Pi f=E_\nu f(X_0)$.
This condition corresponds to the fact that the Markov operator is quasi-compact on the space $\LL$ with $1$ as only
eigenvalue of modulus one and simple (see Hennion and Herv\'e \cite{HenHer01}). Since $\LL\hookrightarrow L^\infty$,
the following result is a special case of Corollary 2 of Durieu \cite{Dur08b}.
\begin{proposition}\label{od}
If $\left(X_n\right)_{n\ge 0}$ is an $\LL$-geometrically ergodic Markov chain then (\ref{eq:4th-moment}) holds
for all $f\in\LL$ such that $Ef(X_0)=0$, with $\alpha=3$ and $\beta=2$.
\end{proposition}

The same is true for dynamical systems whose Perron-Frobenius operators satisfy (\ref{geom}).

\medskip

This gives a large class of examples where our result applies. 

\subsection*{Linear processes}

Let $(A,\|.\|_A)$ be a separable Banach space and $\A$ its Borel sigma algebra. Let $(a_i)_{i\ge0}$ be a sequence of linear forms on $A$ such that there exist
$C>0$ and $0<\theta<1$ such that 
\begin{equation}\label{ai}
|a_i|\le C\theta^i,
\end{equation}
where $|a_i|=\sup_{\|x\|_A\le 1}|a_i(x)|$.
Let $(e_i)_{i\in\Z}$ be an i.i.d. bounded random sequence with values in a compact subset $B\subset A$ and marginal distribution $\mu$. We define the
real-valued linear process $(X_k)_{k\ge0}$ by
$$
X_k=\sum_{i\ge 0}a_i(e_{k-i}),\;k\ge0.
$$
Several results have already been established for empirical processes of linear processes (see Doukhan and Surgailis \cite{DouSur98},
Wu \cite{Wu08}, Dedecker and Prieur \cite{DedPri07}). Here, assumption on the $(a_i)_{i\ge 0}$ is stronger than in the mentioned papers,
but there will be no assumption on the distribution of the $e_i$'s and assumption on the distribution function of
$X_0$ will be weaker.
Note that $(X_k)_{k\ge 0}$ can be viewed as a functional of a Markov chain. 

Let
$Y_k=(e_k,e_{k-1},\dots)$, then $(Y_k)_{k\ge 0}$ is a stationary Markov chain on $B^\N$ (with stationary measure $\mu^{\otimes\N}$) and $X_k=\Phi(Y_k)$ where
$$
\Phi:B^\N\longrightarrow\R,\;\Phi(x_0,x_1,\dots)=\sum_{i\ge 0}a_i(x_i).
$$
Let $Q$ be the Markov transition operator of the chain.
On $B^\N$, we define a metric $d$  by
$$
d(x,y)=\sum_{i\ge 0}\theta^i\|x_i-y_i\|_A
$$
where $x=(x_i)_{i\ge 0}$ and $y=(y_i)_{i\ge0}$. As $B$ is compact, then $(B^\N,d)$ is also compact.
Let us denote by $\LL$ the space of all Lipschitz functions from $B^\N$ to $\R$ provided with the norm $\|.\|$ defined
by
$$
\|f\|=\sup_{x\in B^\N}|f(x)|+\sup_{x\ne y}\frac{|f(x)-f(y)|}{d(x,y)}.
$$
For all $f\in\LL$ and for all $x=(x_i)_{i\ge 0}$ and $y=(y_i)_{i\ge0}\in B^\N$, we have
\begin{eqnarray*}
 |Q^kf(x)-Q^kf(y)|&=&
|E(f(Y_k)|Y_0=x)-E(f(Y_k)|Y_0=y)|\\
&=& |E(f(e_k,\dots,e_1,x_0,\dots))-E(f(e_k,\dots,e_1,y_0,\dots))|\\
&\le &\|f\|E\{d((e_k,\dots,e_1,x_0,\dots),(e_k,\dots,e_1,y_0,\dots))\}\\
&=&C\theta^k\|f\|d(x,y),
\end{eqnarray*}
and
\begin{eqnarray*}
 |Q^kf(x)-Ef(Y_0)|&=&|E(f(Y_k)|Y_0=x)-Ef(Y_k)|\\
&\le&E |f(e_k,e_{k-1},\dots,e_1,x_0,\dots)-f(e_k,e_{k-1},\dots)|\\
&\le&C\theta^k\|f\|E\{d(x,Y_0)\}.
\end{eqnarray*}
Then, we have for all $f\in\LL$,
\begin{eqnarray*}
\|Q^kf-E(f(Y_0))\|&\le& C\theta^k\|f\|.
\end{eqnarray*}
Since $(\LL,\|.\|)\subset (L^\infty(\mu^{\otimes\N}),\|.\|_\infty)$,
by Proposition \ref{od}, $(Y_k)_{k\ge 0}$ satisfies the 4-th moment bound (\ref{eq:4th-moment}) with $\alpha=3$ and $\beta=2$ for all Lipschitz functions.
Further, for all $f\in\LL$ the sequence $\sum_{i=0}^nQ^if(Y_0)$ converges in $L^2(\mu^{\otimes\N})$ and so by Gordin's theorem
(see Gordin \cite{Gor69}), the CLT (\ref{eq:lip-clt}) is satisfied.
Clearly, the function $\Phi$ is a Lipschitz continuous function on $B^\N$, and for all Lipschitz function $g:\R\longrightarrow\R$, $g\circ\Phi$ is also a Lipschitz continuous function on $B^\N$.
Thus conditions (\ref{eq:lip-clt}) and (\ref{eq:4th-moment}) hold for the process $(X_k)_{k\ge0}$, for all Lipschitz function on $\R$. Then Theorem \ref{thm1} applies and we have
\begin{corollary}
Let $(X_k)_{k\ge 0}$ be a real linear process defined by a sequence of linear forms $(a_i)_{i\ge0}$ and a sequence of i.i.d.
bounded random variables $(e_i)_{i\in\Z}$, both on a measurable Banach space $A$. Assume $(a_i)$ satisfies (\ref{ai}) and the distribution function $F$ of $X_0$ satisfies
$$
\omega_F(\delta)\leq D|\log(\delta)|^{-\gamma}\mbox{ for some } D>0 \mbox{ and } \gamma>2.
$$
Then $(U_n(t))_{t\in\R}$ converges in distribution to a mean-zero Gaussian process.
\end{corollary}

In the paper by Dedecker and Prieur \cite{DedPri07}, Corollary 1, $X_0$ has a bounded density. Here, the existence of a density is not needed. Our result is comparable to a result of Wu and Shao \cite{WuSha04}.

\medskip

For a concrete example, consider $A=\{0,1\}$, $a_i=\frac{2}{3^i}$, $i\ge 0$ and $e_k=0$ or $1$ with probability $\frac{1}{2}$, $k\in\Z$.
Then 
$$
X_k=2\sum_{i\ge 0}\frac{e_{k-i}}{3^i},\;k\ge 0
$$
is a stationary process with values in $[0,1]$ and the common distribution function of all the $X_k$ is the Cantor function, which is not absolutly continuous but
$\frac{\log2}{\log3}$-H\"older continuous (see Dovgoshey et al. \cite{DovMarRyaVuo06}).

\subsection*{Expanding maps}

In the setting of expanding maps of the interval, empirical process invariance principles have been established in Collet, Martinez and Schmitt \cite{ColMarSch04} and Dedecker and Prieur \cite{DedPri07} for classes of Lasota-Yorke transformations.
For these maps, the transfer operator has a spectral gap on the space BV of bounded variation functions.  According to Gou\"ezel \cite{Gou08}, there exist some uniformly expanding maps of the interval for which the transfer operator does not act continuously on the space BV, but admits a spectral gap on the space of Lipschitz functions.
The example given by Gou\"ezel is a transformation of the interval $[0,1)$. Let $(a_n)_{n\ge 1}$ be a sequence of positive numbers with $\sum a_n<\frac{1}{4}$ and let $N>0$ be an integer. Denote by $I_n$ the subintervals
$[4\sum_{i=1}^{n-1}a_i,4\sum_{i=1}^{n}a_i)$. We decompose $I_n$ into two subintervals of lenght $2a_n$ denoted by $I_n^{(1)}$
and $I_n^{(2)}$. We can find a map $v_n$ (resp. $w_n$) on $[0,1)$ with image $I_n^{(1)}$ (resp. $I_n^{(2)}$) such that the derivative at a point $x$ is equal to $a_n(1+2\cos^2(2\pi n^4x))$ (resp. $a_n(1+2\sin^2(2\pi n^4x))$). The map $T$ is defined on $I_n$ in such a way that $v_n$ and $w_n$ are two inverse branches of it. It remains the interval
 $[4\sum_{i=1}^{\infty}a_i,1)$ that we subdivide into $N$ subintervals of equal lenght. 
$T$ is defined as an affine transformation on each of these subintervals onto $[0,1)$.

\newtheorem{theo}{Theorem}
\renewcommand{\thetheo}{\empty}

\begin{theo}[Gou\"ezel \cite{Gou08}]
If $a_n=\frac{1}{100n^3}$ and $N=4$, then $T$ is a Lebesgue measure preserving transformation and its associated transfer operator has a spectral gap on the space of Lipschitz functions with a simple eigenvalue at 1 and no other eigenvalue of modulus 1. Further, the transfer operator does not act continuously on BV.
\end{theo}

In this situation, the 4-th moment bound (\ref{eq:4th-moment}) holds and Theorem \ref{thm1} can be used to get an invariance principle for the associated empirical process.

\subsection*{Further applications}

Durieu \cite{Dur08b} has also given 4-th moment bounds for subshifts of finite type, using the Ruelle-Perron-Frobenius theorem,
as in Parry and Pollicott \cite{ParPol90}. Our result thus also applies here.

\medskip

Another application concerns random iterative Lipschitz models and, as a special case, nonlinear autoregressive models $(X_n)_{n\in\N}$ define as follows.
For a real-valued random variable $X_0$ and a given function $f:\R\longrightarrow\R$, let
$$
X_n=f(X_{n-1})+Y_n,\quad n\ge 1,
$$
where and $(Y_n)_{n\ge 1}\subset\R$ is an i.i.d. sequence of $\R$-valued random variables independent of $X_0$.
Such models are studied, e.g., in nonlinear time series analysis.
See Hennion and Herv\'e \cite{HenHer01} Thm.X.16 for conditions under which $(X_n)_{n\in\N}$ is $\LL$-geometrically ergodic.

\paragraph{Acknowledgement} We are grateful to Lo\"{i}c Herv\'e  for 
several lectures introducing us to the spectral gap method, and to J\'er\^ome 
Dedecker for his critical comments on an earlier version of this paper.

\bibliographystyle{plain}

\end{document}